\documentclass[11pt]{article}
\usepackage{amssymb,amsfonts,amsmath,amsthm}
\usepackage{graphicx}
\usepackage[english]{babel}

\newtheorem{theorem}{\indent Theorem}[section]
\newtheorem{lemma}[theorem]{\indent Lemma}

\newtheorem{remark}[theorem]{\indent Remark}
\newtheorem{remarks}[theorem]{\indent Remarks}
\newtheorem{proposition}[theorem]{\indent Proposition}
\newtheorem{corollary}[theorem]{\indent Corollary}
\newtheorem{conjecture}[theorem]{\indent Conjecture}

\newenvironment{definition}{{\bf Definition.}~ }{}

\newcommand{\RR}{\mathbb{R}}

\newcommand{\QQ}{\mathbb{Q}}
\newcommand{\KK}{\mathbb{K}}

\newcommand{\ZZ}{\mathbb{Z}}

\newcommand{\HH}{\mathbb{H}}

\newcommand{\FF}{\mathcal{F}}

\newcommand{\SL}{\mathrm{SL}}
\newcommand{\PSL}{\mathrm{PSL}}
\newcommand{\PGL}{\mathrm{PGL}}

\newcommand{\PSO}{\mathrm{PSO}}

\newcommand{\Id}{\mathrm{Id}}
\newcommand{\tr}{\mathrm{tr}}

\newcommand{\im}{\mathrm{Im}}
\newcommand{\GA}{\Gamma}
\newcommand{\ga}{\gamma}

\newcommand{\ol}{\overline}
\newcommand{\dH}{\partial\mathbb{H}}

\begin{document}

\title{Weyl chamber flow on irreducible quotients of products of $\PSL(2,\RR)$}
\author{Damien \textsc{Ferte}\\
}
\date{}
\maketitle
\noindent{\bf Abstract: }We study the topological dynamics of the action of the diagonal subgroup on quotients $\GA\backslash\PSL(2,\RR)\times\PSL(2,\RR)$, where $\GA$ is an irreducible lattice. Closed orbits are described and a set of points of dense orbit is explicitly given. Such properties are expressed using the Furstenberg boundary of the symmetric space $\HH\times\HH$.\\

\vskip 10mm

\section{Introduction}
Let $\rm G$ be a real connected semi-simple Lie group, $\rm K$ a maximal compact subgroup, $\rm A$ a Cartan subgroup and $\rm M$ the centralizer of $\rm A$ in $\rm K$. If $\GA$ is a lattice of $\rm G$, the action (by right-translation) of $\rm A$ on the coset $\GA\backslash{\rm G}/{\rm M}$ is the {\it Weyl chamber flow} over the locally symmetric space $\GA\backslash{\rm G}/{\rm K}$. When the group ${\rm G}$ is $\PSL(2,\RR)$, this action is conjugated to the geodesic flow on the unit tangent bundle of the finite-volume hyperbolic surface $\GA\backslash\HH$. In this case, there are a lot of different kinds of closed $\rm A$-invariant subsets. By contrast, when the dimension of the Cartan subgroup $\rm A$ (the rank of $\rm G$) is greater than 2, the situation is expected to be very rigid. In \cite{Margulis survey1}, G.\,A.\,Margulis conjectured that, excluding a situation of factorisation by a rank one action, each $\rm A$-orbit closure is {\it algebraic} i.e. is the orbit of a closed connected subgroup of $\rm G$ containing $\rm A$. Points of $\GA\backslash\rm G$ whose orbit is closed or dense satisfy this conjecture. Closed and compact orbits come from unipotent or semi-simple elements of the lattice $\GA$. (See \cite{Tomanov-Weiss} for a description of closed orbits and \cite{Prasad-Raghunathan} for a condition of compactness.) Moore's ergodicity theorem (see \cite{Zimmer}) implies that, for the finite Haar measure, almost all point in $\GA\backslash\rm G$ have a dense orbit. But there is no explicit way to find or construct such dense orbits. For the groups $\SL(n,\RR)$ and their products, E. Lindenstrauss and B. Weiss proved that each point whose $\rm A$-orbit closure contains a compact orbit also satisfies Margulis' conjecture. See also \cite{Mozes} for the study of the group $\PGL(2,\QQ_p)\times\PGL(2,\QQ_l)$.

This article deals with the special case of the group ${\rm G}=\PSL(2,\RR)\times\PSL(2,\RR)$. The subgroup $\rm A$ is the maximal diagonal subgroup and we will also consider the semi-subgroup
$${\rm A^+}=\left\{\left(\pm\begin{pmatrix}\lambda_1&0\\0&\lambda_1^{-1}\end{pmatrix}, \pm\begin{pmatrix}\lambda_2&0\\0&\lambda_2^{-1}\end{pmatrix}\right): \lambda_1,\lambda_2\geq 1\right\}.$$
If the lattice $\GA$ in $\rm G$ is {\it reducible} (up to finite index, it is a product of two lattices of $\PSL(2,\RR)$), this is the simplest case of factorisation: every closed invariant subset for the action of $\rm A$ is the product of closed invariant subsets for the geodesic flows. Therefore we shall assume that $\GA$ is an {\it irreducible} lattice of $\rm G$. For instance, if $\KK$ is a real quadratic field (of Galois automorphism $\sigma$) and $\cal O$ is its integers ring, the injection of the group $\PSL(2,\mathcal O)$ in $\rm G$ given by
$$\pm\begin{pmatrix}a&b\\c&d\end{pmatrix}\longmapsto \left(\pm\begin{pmatrix}a&b\\c&d\end{pmatrix}, \pm\begin{pmatrix}a^\sigma&b^\sigma\\c^\sigma&d^\sigma\end{pmatrix}\right)$$
is an irreducible lattice called {\it Hilbert modular lattice} associated to $\KK$. For irreducible lattices of $\rm G$, Margulis' conjecture can be strentghtened in the following way.

\begin{conjecture}[Margulis]\label{margulis conjecture}Let $\GA$ be an irreducible lattice of $\rm G=\PSL(2,\RR)\times\PSL(2,\RR)$ and $\rm A$ be the maximal diagonal subgroup of $\rm G$. Let $x$ be a point in $\GA\backslash\rm G$. Thus the orbit $x\rm A$ is either closed or dense in $\GA\backslash\rm G$.
\end{conjecture}

In \cite{Lindenstrauss}, the question of classification of $\rm A$-ergodic finite measures on $\GA\backslash\rm G$ is studied and related to a ``quantum unique ergodicity'' conjecture.

In this article, we give an explicit set of points whose orbit (by the semi-group $\rm A^+$) is dense. The condition is expressed on the factors $\dH=\RR\cup\{\infty\}$ of the Furstenberg boundary $\cal F=\dH\times\dH$ of the symmetric space $\HH\times\HH$ associated to the group $\rm G$.

We use here the terminology of \cite{Shimizu}. An element of $\rm G$ is said to be {\it hyperbolic} (resp. {\it parabolic}, {\it elliptic}) if both components are hyperbolic (resp. parabolic, elliptic) elements of $\PSL(2,\RR)$ (we use the convention that the unit element is neither hyperbolic, nor parabolic, nor elliptic). A non-trivial element is said to be {\it mixed} if its two components are not of the same kind. A hyperbolic element $\ga=(\ga_1,\ga_2)$ is said to be {\it hyper-regular} if the dominant (positive) eigenvalues of $\ga_1$ and $\ga_2$ are distinct (in this situation, this definition coincides with the definition of ``hyper-regular'' in \cite{Prasad-Raghunathan}). An example of non hyper-regular elements is given by the canonical injection of the subgroup $\PSL(2,\ZZ)$ in a Hilbert modular lattice. The group $\rm G$ acts naturally on $\cal F$, but also on each factor $\dH$ by the corresponding component. For instance, a point $\xi_1$ in the first factor $\dH$ of $\cal F=\dH\times\dH$ is fixed by an element $\ga=(\ga_1,\ga_2)$ if $\ga_1(\xi_1)=\xi_1$.

\begin{theorem}\label{sufficient conditions of density} Let $\GA$ be an irreducible lattice of $\rm G=\PSL(2,\RR)\times\PSL(2,\RR)$. Let $g=(g_1,g_2)$ belong to the group $\rm G$ and $x$ be the class of $g$ in $\GA\backslash\rm G$. Then the semi-orbit $x\rm A^+$ is dense in $\GA\backslash\rm G$ if one of the following conditions holds:\\
\indent 1) a mixed element of $\GA$ fixes one of the points $g_1(\infty)$ and $g_2(\infty)$ of $\dH$,\\
\indent 2) a hyperbolic hyper-regular element of $\GA$ fixes exactly one of the points $g_1(\infty)$ and $g_2(\infty)$.
\end{theorem}

We also describe compact orbits and retrieve conditions of compactness of \cite{Prasad-Raghunathan}.

\begin{theorem}\label{compact orbits} Let $\GA$ be an irreducible lattice of $\rm G=\PSL(2,\RR)\times\PSL(2,\RR)$. Let $g$ belong to the group $\rm G$ and $x$ be the class of $g$ in $\GA\backslash\rm G$. Then the following properties are equivalent:\\
\indent 1) the orbit $x\rm A$ is compact,\\
\indent 2) the subgroup $g^{-1}\GA g\cap{\rm A}$ is isomorphic to $\ZZ^2$,\\
\indent 3) the points $g(\infty,\infty)$ and $g(0,0)$ of $\FF$ are fixed by a hyperbolic hyper-regular element of $\GA$.
\end{theorem}

When the lattice is not uniform, there exist closed, but non-compact orbits, coming from the parabolic points on $\cal F$ (see paragraph \ref{Parabolic points and horoballs} for definition of conjugate parabolic points). In \cite{Tomanov-Weiss}, a general statement is given in algebraic terms.

\begin{theorem}\label{closed non-compact orbits} Let $\GA$ be an irreducible (non-uniform) lattice of $\rm G=\PSL(2,\RR)\times\PSL(2,\RR)$. Let $g$ belong to the group $\rm G$ and $x$ be the class of $g$ in $\GA\backslash\rm G$. Then the following properties are equivalent:\\
\indent 1) the orbit $x\rm A$ is closed and non-compact,\\
\indent 2) the subgroup $g^{-1}\GA g\cap{\rm A}$ is isomorphic to $\ZZ$,\\
\indent 3) the points $g(\infty,\infty)$ and $g(0,0)$ (or $g(0,\infty)$ and $g(0,\infty)$) are conjugate parabolic points of $\FF$, relatively to $\GA$.
\end{theorem}

\begin{remarks}\label{remarks in introduction} 1) The element $\footnotesize{\pm\begin{pmatrix}0&-1\\1&0\end{pmatrix}}$ of the group $\PSL(2,\RR)$ normalizes the diagonal subgroup of $\PSL(2,\RR)$ and exchanges the points $0$ and $\infty$. Then, the previous theorem \ref{sufficient conditions of density} give sufficient conditions of density for the full diagonal group $\rm A$ using the points $g_1(0),g_1(\infty),g_2(0)$ and $g_2(\infty)$ of $\dH$.\\
2) If a hyperbolic hyper-regular element of $\GA$ fixes both points $g_1(\infty)$ and $g_2(\infty)$ then the semi-orbit $x\rm A^+$ is ``asymptotic'' to a compact orbit and hence cannot be dense.
\end{remarks}

We obtain the following corollary for uniform lattices.

\begin{corollary}\label{corollary for uniform lattices} Let $\GA$ be a uniform irreducible lattice of $\rm G=\PSL(2,\RR)\times\PSL(2,\RR)$. Then conjecture \ref{margulis conjecture} is true for every point $x=\GA g$ in $\GA\backslash\rm G$ such that one of the points $g_1(0),g_1(\infty),g_2(0),g_2(\infty)$ of $\dH$ is fixed by a non-trivial element of $\GA$.
\end{corollary}

This article is organized as follows. Paragraph 1 contains the useful results about irreducible lattices. After some recalls about the geometry of the locally symmetric space $\GA\backslash\HH\times\HH$ (paragraph 2), theorems \ref{compact orbits} and \ref{closed non-compact orbits} are proved in paragraph 3. Theorem \ref{sufficient conditions of density} and corollary \ref{corollary for uniform lattices} are proved in paragraph 4. We use classical facts about the geometry of the hyperbolic plane: geodesics, compactification, dynamics of isometries on $\HH$ and on the boundary $\dH$. If all notations coincide with the upper half-plane model of the hyperbolic plane, we use the disc model for the figures. 

\section{Properties of irreducible lattices}

Denote by $p_i$, $i=1,2,$ the projections
$$p_i: \mathrm G\longrightarrow\PSL(2,\RR),\,(g_1,g_2)\longmapsto g_i$$
on each factor. If $\GA$ is an irreducible lattice of $\rm G$, the subgroups $p_i(\GA)$ are dense in $\PSL(2,\RR)$ and the kernels $\GA\cap\ker p_i$ are central in $\ker p_i$ (which is isomorphic to $\PSL(2,\RR)$) and therefore
 trivial (see \cite{Raghunathan}).\\
\indent Concerning the elements of irreducible lattices, we have the following. It is well-known that hyperbolic (and hyper-regular) elements exist in such a lattice. The density of the projections implies that an irreducible lattice $\GA$ contains elements with an elliptic component. If $\GA$ is a uniform lattice of $\rm G$, then every element of $\GA$ is semi-simple. Thus $\GA$ doesn't contain elements with a parabolic component. If $\ga=(\ga_1,\ga_2)$ belongs to a Hilbert modular lattice, then we have $\tr(\ga_2)=\tr(\ga_1)^\sigma$. Hence $\ga_1$ is parabolic if and only if $\ga_2$ is also parabolic. Therefore, if $\GA$ is a uniform irreducible lattice or a Hilbert lattice, the non-trivial elements of $\GA$ are of the following kinds:\\
- hyperbolic (and there exist hyper-regular elements),\\
- elliptic of finite order (if and only if $\GA$ is not torsion free),\\
- mixed (one component is hyperbolic and the other is elliptic),\\
- parabolic (if $\GA$ is not uniform).\\
For general non-uniform lattices, A. Selberg proved in \cite{Selberg} his arithmeticity theorem:

\begin{theorem}[Selberg]\label{Selberg arithmeticity theorem} A non-uniform irreducible lattice $\GA$ of $\rm G$ is, up to conjugation by an element of $\PGL(2,\RR)\times\PGL(2,\RR)$, commensurable to a Hilbert modular lattice.
\end{theorem}

Hence the previous classification of elements is valid for any irreducible lattice of $\rm G$.

The following result can be seen like a refinement of the density of the projections of an irreducible lattice on each factor of $\rm G$. This statement is not symmetric but also true with a permutation of indices.

\begin{proposition}\label{refined density} Let $\GA$ be an irreducible lattice of $\rm G$. Let $g_1$ belong to $\PSL(2,\RR)$ and $\eta_2^+,\eta_2^-$ belong to $\dH$. Then there exists a sequence $(\ga_n)_n$ in the lattice $\GA$ satisfying:\\
1) $\displaystyle{\lim_{n\rightarrow+\infty}p_1(\ga_n)=g_1}$,\\
2) for every point $z_2$ in $\HH$: $\displaystyle{\lim_{n\rightarrow+\infty} p_2(\ga_n)z_2=\eta_2^+}$ and $\displaystyle{\lim_{n\rightarrow+\infty} p_2(\ga_n)^{-1}z_2=\eta_2^-}$,\\
3) for every point $\xi_2$ of $\dH$, distinct from $\eta_2^-$: $\displaystyle{\lim_{n\rightarrow+\infty} p_2(\ga_n)(\xi_2)=\eta_2^+}$.
\end{proposition}

\begin{proof} Recall that it is sufficient to prove assertion \textit{2)} for one point $z_2$. Let $V_1$ be a neighborhood of $g_1$ in $\PSL(2,\RR)$ and $V_2^+$ (resp. $V_2^-$) be a neighborhood of $\eta_2^+$ (resp. $\eta_2^-$) in $\dH$. There exists an element $\alpha$ in $\GA$ such that $p_1(\alpha)$ is elliptic and $p_2(\alpha)$ is hyperbolic. Let $\xi^+_2$ (resp. $\xi^-_2$) be the attractive (resp. repulsive) fixed point of $p_2(\alpha)$. If $(\alpha_n)_n$ is a appropriate subsequence of the positive powers of $\alpha$, the sequence $(p_1(\alpha_n))_n$ converges to the point $\Id$ of $\PSL(2,\RR)$ and the sequence $(p_2(\alpha_n)z_2)_n$ (resp. $(p_2(\alpha_n)^{-1}z_2)_n$) of points of $\HH$ converges to the point $\xi_2^+$ (resp. $\xi_2^-$) of $\dH$. The projection $p_1(\GA)$ is dense in $\PSL(2,\RR)$, hence there exists $\gamma$ in $\GA$ such that $p_1(\gamma)$ belongs to $V_1$. The points $\xi_2^+$ and $\xi_2^-$ are distinct, the group $\PSL(2,\RR)$ acts transitively on the set of distinct points of $\dH$ and the projection $p_2(\GA)$ is dense in $\PSL(2,\RR)$. Therefore, there exists an element $\delta$ in $\GA$ such that
$$p_2(\delta)\xi_2^+\in V_2^+\quad\text{and}\quad p_2(\delta)\xi_2^-\in p_2(\gamma)V_2^-.$$
The sequence $(\delta\alpha_n\delta^{-1}\gamma)_n$ of elements of $\GA$ satisfies
$$\lim_{n\rightarrow+\infty}p_1(\delta\alpha_n\delta^{-1}\gamma)=p_1(\gamma)\in V_1,$$
$$\lim_{n\rightarrow+\infty}p_2(\delta\alpha_n\delta^{-1}\gamma)z_2=p_2(\delta)(\xi_2^+)\in V_2^+,$$
$$\lim_{n\rightarrow+\infty} p_2(\delta\alpha_n\delta^{-1}\gamma)^{-1}z_2= p_2(\gamma)^{-1}p_2(\delta)(\xi_2^-)\in V_2^-.$$
This proves the existence of a sequence $(\ga_n)_n$ in the lattice $\GA$ satisfying assertions \textit{1)} and \textit{2)}.\\
It remains to prove assertion {\it 3)}. This is a consequence of {\it 2)}. To prove this, we can omit the index ``2'' which is not useful. Consider the set $\cal E$ of points $\xi$ in $\dH$ such that the sequence $(\ga_n(\xi))_n$ does not converge to $\eta^+$. If this set contains two distincts points, they can be joined by a geodesic line and any point $z$ on this geodesic cannot satisfies
$$\lim_{n\rightarrow+\infty}\ga_nz=\eta^+.$$
Hence $\cal E$ contains at most one point $\xi$ and it remains to show that this point (if it exists) is $\eta^-$. Assume this is false and let $V$ be a (small) neighborhood of $\eta^+$ in $\dH$. Then, for $n$ large enough, $\ga_n(\eta^-)$ belongs to $V$, and there exists a neighborhood $W$ of $\eta^-$ in $\dH$ such that $\ga_nW$ is contained in $V$, for all $n$. The points $\xi$ and $\eta^-$ can be joined by a geodesic line. Let $z$ be a point of this line and for any $n$, let $\sigma_n$  be the oriented geodesic line passing through $\ga_n^{-1}z$ and $z$ (in this order). Denote by $\sigma_n(-\infty)$ and $\sigma_n(+\infty)$ the extremities of this geodesic. The points $\ga_n(\sigma_n(-\infty))$ and $\ga_n(\sigma_n(+\infty))$ are the extremities of the geodesic line $\ga_n\sigma_n$ passing through $z$ and $\ga_nz$ (in this order). The point $\ga_n(\sigma_n(-\infty))$ belongs to $V$, and the fact that $(\ga_nz)_n$ converges to $\eta^+$ implies that $\ga_n(\sigma_n(+\infty))$ also belongs to $V$. This is impossible because the geodesic line $\ga_n\sigma_n$ contains the point $z$.
\end{proof}

We obtain the following corollary which allows to simplify the setting of Margulis' conjecture and the proof of Lindenstrauss-Weiss' result in this situation (see \cite{Ferte} for details).

\begin{corollary} Let $\GA$ be an irreducible lattice of $\rm G=\PSL(2,\RR)\times\PSL(2,\RR)$ and $\rm F$ be a closed connected subgroup of $\rm G$ strictly containing the group $\rm A$. Then, in $\GA\backslash\rm G$, every orbit of $\rm F$ is dense.
\end{corollary}

\begin{proof} Let $\rm U^+$ (resp. $\rm U^-$) be the upper (resp. lower) unipotent subgroup of $\rm G$ and ${\rm A}_1\times{\rm A}_2$, ${\rm U}_1^+\times{\rm U}_2^+$ and ${\rm U}_1^-\times{\rm U}_2^-$ be the canonical decomposition of $\rm A$, $\rm U^+$ and $\rm U^-$ in $\rm G=\PSL(2,\RR)\times\PSL(2,\RR)$. If $\rm F$ is a connected subgroup of $\rm G$ strictly containing the group $\rm A$, then its Lie algebra is invariant under the adjoint action of $\rm A$. Thus $\rm F$ contains one of the four unipotent triangular subgroups ${\rm U}_1^+,{\rm U}_1^-,{\rm U}_2^+,{\rm U}_2^-$ of $\rm G$. It is therefore sufficient (up to transposition and permutation of indices) to prove that the subgroup ${\rm A}_2{\rm U}_2^+$ acts minimally on $\GA\backslash\rm G$. We consider the ``dual'' action of $\GA$ on ${\rm G}/{\rm A}_2{\rm U}_2^+$. This coset can be identified with $\PSL(2,\RR)\times\dH$, on which the action of $\GA$ is given by
$$\ga(h_1,\xi_2)=(p_1(\ga)h_1,p_2(\ga)\xi_2).$$
The minimality of the action is therefore an easy consequence of the previous proposition (assertions \textit{1} and \textit{3}).
\end{proof}

\section{Parabolic points and horoballs}\label{Parabolic points and horoballs}

Let $\beta_.(.,.)$ be the classical Busemann cocycle on $\HH$ whose sign is fixed by the following equality
$$\beta_\infty(x,y)=\ln\left(\frac{\im x}{\im y}\right)\qquad\text{for $x,y$ in }\HH.$$
We define on $\HH\times\HH$ the cocycle $\beta_\xi(.,.)$ as the sum of the Busemann cocycles on each factor: if $\xi=(\xi_1,\xi_2)$ belongs to $\cal F=\dH\times\dH$ and $z=(z_1,z_2),z'=(z'_1,z'_2)$ belong to $\HH\times\HH$,
$$\beta_\xi(z,z')=\beta_{\xi_1}(z_1,z'_1)+\beta_{\xi_2}(z_2,z'_2).$$
When the point $\xi$ is $(\infty,\infty)$, the cocycle is then given by
$$\beta_\infty(z,z')=\ln\left(\frac{\im z_1\im z_2}{\im z_1'\im z_2'}\right).$$

\begin{definition}\label{horoball} Fix a point $z_o$ in $\HH\times\HH$. Let $\xi$ be a point in the boundary $\cal F=\dH\times\dH$ and $T$ be a real. The {\it horoball} based at $\xi$ and of level $T$ is the subset
$${\rm HB}(\xi,T)=\{z\in\HH\times\HH\;:\;\beta_\xi(z,z_o)>T\}.$$
\end{definition}

The following lemma will be used in the proof of theorem \ref{closed non-compact orbits}.

\begin{lemma}\label{basepoint of horoball} Fix a point $z=(z_1,z_2)$ in $\HH\times\HH$. Let $\xi=(\xi_1,\xi_2)$ be in $\cal F$, $g=(g_1,g_2)$ be in $\rm G$ and $a=(a_1,a_2)$ be a non-trivial element of $\rm A$. Assume that a horoball ${\rm HB}(\xi,T)$ based at $\xi$ contains infinitely many points $ga^nz$, where $n$ is a positive integer. Then
$$\xi_1=g_1a_1^+\;\text{if}\;a_2=\Id\qquad\text{(resp.\;}\xi_2=g_2a_2^+\;\text{if}\; a_1=\Id\text{)},$$ where $a_i^+$ (which equals $0$ or $\infty$) is the attractive point in $\dH$ of the hyperbolic isometry $a_i$.
\end{lemma}

\begin{proof} It is sufficient to prove the first assertion. If $a_2=\Id$, $a_1$ is non-trivial and the sequence $(a_1^nz_1)_n$ of $\HH$ goes to $a_1^+$ when $n$ goes to infinity. Hence $(\beta_{g_1^{-1}\xi_1}(a_1^nz_1,g_1^{-1}z_1))_n$ goes to $-\infty$ if $g_1^{-1}\xi_1$ is different from $a_1^+$. The point $ga^nz$ belongs to ${\rm HB}(\xi,T)$, therefore we have
$$\beta_{g_1^{-1}\xi_1}(a_1^nz_1,g_1^{-1}z_1)+\beta_{g_2^{-1}\xi_2}(z_2,g_2^{-1}z_2) =\beta_{g^{-1}\xi}(a^nz,g^{-1}z)=\beta_\xi(ga^nz,z)>T.$$
Consequently the point $g_1^{-1}\xi_1$ is equal to $a_1^+$.
\end{proof}

\begin{definition}\label{parabolic point} A point $\xi$ of the boundary $\cal F$ is said to be {\it parabolic} (with respect to a subgroup $\GA$ of $\rm G$) if it is fixed by a parabolic element of $\GA$.
\end{definition}\\

For instance, the point $(\infty,\infty)$ is parabolic with respect to any Hilbert modular lattice of $\rm G$ and its stabilizer in such a lattice is given by the set of matrices
$$\left\{\pm\begin{pmatrix}a&b\\0&a^{-1}\end{pmatrix}\,:\,a\in\mathcal O^\times,b\in\mathcal O\right\}.$$
The Galois conjugate $a^\sigma$ of an element $a$ of the unit group $\mathcal O^\times$ is equal to $\pm a^{-1}$ because $\vert aa^\sigma\vert=1$ hence every hyperbolic element of the stabilizer is conjugate to an element of the subgroup
$$\mathrm A'=\left\{\left(\pm\begin{pmatrix}\lambda&0\\0&\lambda^{-1}\end{pmatrix}, \pm\begin{pmatrix}\lambda^{-1}&0\\0&\lambda\end{pmatrix}\right) \,:\,\lambda\in\mathcal\RR^*\right\}$$
which preserves each horoball based on $(\infty,\infty)$ and all of whose non-trivial element is hyperbolic but non hyper-regular. The following results on stabilizers of parabolic points are proved in \cite{Shimizu}.

\begin{proposition}\label{parabolic stabilizer} Let $\GA$ be an irreducible lattice of $\rm G$, $\xi$ be a parabolic point in $\cal F$ and $\GA_\xi$ the stabilizer of $\xi$ in $\GA$.\\
\indent - There exists a real number $T$ such that:
$$\ga{\rm HB}(\xi,T)\cap{\rm HB}(\xi,T)=\emptyset\quad\text{for every element $\ga$ in $\GA-\GA_\xi$}.$$
\indent - If $\xi=g(\infty,\infty)$, the stabilizer $\GA_\xi$ is a (uniform) lattice in the solvable subgroup $g{\rm A'U^+}g^{-1}$. In particular this subgroup does not contain hyper-regular element and globally preserves each horoball based on $\xi$.
\end{proposition}

Hence a parabolic point is fixed by a hyperbolic element. Two parabolic points are {\it conjugate} if they are fixed by the same hyperbolic element.\\

The following theorem (called {\it Property (F)} in \cite{Shimizu}) which is obvious for uniform lattices and known for Hilbert modular lattices is true for any irreducible lattice of $\rm G$, using the arithmeticity theorem \ref{Selberg arithmeticity theorem}. 

\begin{theorem}\label{fundamental domain} Let $\GA$ an irreducible lattice in $\rm G$. Then there exists a real number $T$ and some representatives $\xi^1,\ldots,\xi^r$ in $\cal F$ of all equivalence classes of parabolic points such that $\GA$ has a fundamental domain $\cal D$ in $\HH\times\HH$ which is a disjoint union
$$\mathcal D=K\sqcup V_1\sqcup\ldots\sqcup V_r,$$
where $K$ is a compact subset and $V_i$ is a fundamental domain for the action of $\GA_{\xi^i}$ on the horoball ${\rm HB}(\xi^i,T)$.
\end{theorem}

\begin{corollary}\label{disjoint union of horospheres} The complement of $\GA K$ in $\HH\times\HH$ is a disjoint union of horoballs of level $T$ based on all parabolic points of $\cal F$.
\end{corollary}

\section{Closed orbits}

In this section, we prove theorems \ref{closed non-compact orbits} and \ref{compact orbits} using the following classical result. Let $\GA$ be a discrete subgroup of $\rm G$ and $\rm A$ be the maximal diagonal subgroup of $\rm G$. Let $g$ be an element of $\rm G$ and $x$ the class of $g$ in $\GA\backslash\rm G$. Then the (well-defined) map
$$\Psi_x:(g^{-1}\GA g\cap{\rm A})\backslash{\rm A}\longrightarrow\GA\backslash{\rm G},\,(g^{-1}\GA g\cap{\rm A})a\longmapsto xa$$
is continuous, injective and its image is precisely the orbit $x\rm A$.

\begin{lemma}\label{properness} Assume moreover that $\GA$ is an irreducible lattice of $\rm G$. With these notations:\\
1) The orbit $x\rm A$ is closed in $\GA\backslash\rm G$ if and only if the application $\Psi_x$ is proper,\\
2) If the subgroup $g^{-1}{\GA}g\cap\rm A$ is non-trivial, 1) holds.
\end{lemma}

\begin{proof} \textit{1)} If $\Psi_x$ is proper, then $x\rm A$ is closed. Conversely, assume that $x\rm A$ is closed and $\Psi_x$ is not proper: there exists a sequence $(a_n)_n$ in $\rm A$ diverging in $(g^{-1}\GA g\cap{\rm A})\backslash{\rm A}$ but such that the sequence $(xa_n)_n$ converges in $\GA\backslash\rm G$ to a point $xa$, with $a$ in $\rm A$. Then there exists a sequence $(\ga_n)_n$ in $\GA$ such that $(\ga_nxa_n)_n$ converges in $\rm G$ to $ga$, hence in $\rm G/A$ to the class $g\rm A$ of $g$ in $\rm G/A$. But the set $\GA g\rm A$ is assumed to be closed in $\rm G$, therefore the $\GA$-orbit of the class $g\rm A$ is closed hence discrete (a consequence of the countability of $\GA$ and the Baire's property of the quotient set $\rm G/\rm A$). So the sequence $(\ga_ng{\rm A})_n$ is stationary in $\rm G/\rm A$: there exists $a'$ in $\rm A$ and indices $m>n$ such that
$$\ga_m\neq\ga_n\quad\text{and}\quad\ga_mg=\ga_nga'.$$
Hence $g^{-1}\GA g\cap{\rm A}$ is non-trivial and it remains to prove the second part of the lemma.\\
\textit{2)} Assume there exists a non-trivial element $\ga$ in $\GA$ such that $g^{-1}\ga g$ belongs to $\rm A$. By the hypothesis of irreducibility, the centralizer of $g^{-1}\ga g$ in $\rm G$ is precisely $\rm A$. The properness of $\Psi_x$ follows from the following. Assume that the sequence $(ga_n)_n$ is convergent modulo $\GA$: there exists a point $h$ in $\rm G$ and a sequence $(\ga_n)_n$ in $\GA$ such that $(\ga_nga_n)_n$ converges to $h$. The sequence of elements
$$\ga_n\ga\ga_n^{-1}=\ga_nga_n(g^{-1}\ga g)a_n^{-1}g^{-1}\ga_n^{-1}$$
converges to $hg^{-1}\ga gh^{-1}$. Thus the sequence $(\ga_n\ga\ga_n^{-1})_n$ is stationary: there exists a index $n_o$ such that
$$\ga_n\ga\ga_n^{-1}=\ga_{n_o}\ga\ga_{n_o}^{-1}.$$
Consequently, each element $\ga_{n_o}^{-1}\ga_n$ commutes with $\ga$. Therefore $g^{-1}\ga_{n_o}^{-1}\ga_ng$ belongs to $g^{-1}\GA g\cap\rm A$. The sequence $(g^{-1}\ga_{n_o}^{-1}\ga_nga_n)_n$ converges to $g^{-1}\ga_{n_o}^{-1}h$ thus $(a_n)_n$ is convergent modulo $g^{-1}\GA g\cap\rm A$
\end{proof}

\begin{proof}[Proof of theorem \ref{closed non-compact orbits}]
We denote by $\Lambda$ the discrete subgroup $g^{-1}\GA g\cap{\rm A}$ of $\rm A$.\\
${\it 1)}\Rightarrow{\it 3)}$ We assume that $x\rm A$ is closed but non-compact. Then $\Lambda$ is not a uniform lattice in $\rm A$. Let $z_o$ be a point in $ \HH\times\HH$ and consider the following proper map:
$$\pi:\GA\backslash{\rm G}\rightarrow\GA\backslash\HH\times\HH,\,\GA g\longmapsto\GA gz_0.$$
The orbit $x\rm A$ is closed and non-compact, thus by lemma \ref{properness} the map $\Psi_x$ is proper and the subset $(\pi\circ\Psi_x)^{-1}(K)$ is a compact subset of the non-compact set $\Lambda\backslash\rm A$.\\
\indent Assume $\Lambda$ is trivial. Then $(\pi\circ\Psi_x)^{-1}(K)$ is a compact subset of $\rm A$. Let $C$ be the unique unbounded connected component of ${\rm A}-(\pi\circ\Psi_x)^{-1}(K)$; thus the subset $gCz_o$ of $\HH\times\HH$ is connected and disjoint from $\GA K$. Therefore, by corollary \ref{disjoint union of horospheres}, it is contained in a horoball ${\rm HB}(\xi,T)$ based at a parabolic point. There exists an element $a=(a_1,\Id)$ in $C$, $a_1\neq\Id$, such that the elements $a^n$ and $a^{-n}$ belong to $C$ for every positive integer $n$. The points $ga^nz_o$ and $ga^{-n}z_o$ therefore belong to ${\rm HB}(\xi,T)$. The lemma \ref{basepoint of horoball} imply that the point $\xi_1$ is equal to $g_1(\infty)$ and $g_1(0)$. This is a contradiction. Consequently the discrete subgroup $\Lambda$ of $\rm A$ is not trivial. It is therefore isomorphic to $\ZZ$ because $\Lambda\backslash\rm A$ is not compact.\\
\indent The set $(\pi\circ\Psi_x)^{-1}(K)$ is a compact subset of the cylinder $\Lambda\backslash\rm A$, thus it is contained in a subset $\Lambda\backslash B$ of $\Lambda\backslash\rm A$ where $B$ is a ``band'' in $A$, invariant by the subgroup $\Lambda$. The connected components $C^+$ and $C^-$ of ${\rm A}-B$ satisfy
$$gC^+z_o\subseteq{\rm HB}(\xi^+,T)\quad\text{and}\quad gC^-z_o\subseteq{\rm HB}(\xi^-,T)$$
for some parabolic points $\xi^+$ and $\xi^-$. The group $\Lambda$ is generated by an element 
$g^{-1}\ga g$ where $\ga$ belongs to $\GA$. Its components are non-trivial (by irreducibility of $\GA$), therefore $C^+$ and $C^-$ contain non-trivial elements of the form $(a_1,\Id)$ and $(\Id, a_2)$. Applying once again the lemma \ref{basepoint of horoball} with such elements, we obtain
$$\{\xi^+,\xi^-\}=\{g(0,0),g(\infty,\infty)\}\quad\text{or}\quad\{g(0,\infty),g(0,\infty)\},$$
according to the position of $C^+$ and $C^-$. Moreover, the points $\xi^+$ and $\xi^-$ are fixed by $\ga=gag^{-1}$.\\
${\it 3)}\Rightarrow{\it 2)}$ If an element $\ga$ of $\GA$ fixes the parabolic points $g(0,0)$ and $g(\infty,\infty)$ (or $g(0,\infty)$ and $g(0,\infty)$), the diagonalizable subgroup $\GA\cap g{\rm A}g^{-1}$ is non-trivial and it cannot be isomorphic to $\ZZ^2$ since it is contained in the stabilizer of a parabolic point (proposition \ref{parabolic stabilizer}).\\
${\it 2)}\Rightarrow{\it 1)}$ Under the assumption {\it 2)}, the map $\Psi_x$ is proper by lemma \ref{properness}. Thus $x{\rm A}$ is closed, but not compact since $\Lambda$ is not a lattice in $\rm A$.
\end{proof}

\begin{proof}[Proof of the theorem \ref{compact orbits}] Recall that $\Lambda$ denote the subgroup $g^{-1}\GA g\cap\rm A$.\\
${\it 1)}\Rightarrow{\it 2)}$ The orbit is compact, therefore the map $\Psi_x$ is proper and $\Lambda\backslash\rm A$ is compact. The discrete subgroup $\Lambda$ has to be a lattice in $\rm A$.\\
${\it 2)}\Rightarrow{\it 3)}$ Since any lattice in $\rm A$ contains a hyper-regular hyperbolic element, the lattice $\Lambda$ contains an element $g^{-1}\ga g$ where $\ga$ is a hyper-regular element of $\GA$. Thus $\ga$ fixes both points $g(\infty,\infty)$ and $g(0,0)$.\\
${\it 3)}\Rightarrow{\it 1)}$ Let $\ga$ be a hyper-regular element fixing $g(\infty,\infty)$ and $g(0,0)$. Then $g^{-1}\ga g$ belongs to $\Lambda$ and the orbit $x\rm A$ is closed by lemma \ref{properness}. This orbit cannot be non-compact by theorem \ref{closed non-compact orbits} since $\ga$ is not contained in the stabilizer of a parabolic point.
\end{proof}

\section{Dense orbits}

In this section, we prove theorem \ref{sufficient conditions of density}. We use here the ``geometric language'': an element of $\PSL(2,\RR)$ can be seen as a unit tangent vector (or a geodesic ray) of the hyperbolic plane. Therefore an element of $\rm G=\PSL(2,\RR)\times\PSL(2,\RR)$ can be seen as a product of two unit tangent vectors of $\HH$ (or a Weyl chamber in $\HH\times\HH$). The action (by right translation) of the diagonal group $\rm A$ (resp. the diagonal semigroup $\rm A^+$) on $\rm G$ is then conjugated to the product of the geodesic flows (resp. in positive time). We will use the notation $\phi^t$ for the geodesic flow at time $t$. This is in fact the diagonal matrix $$\pm\begin{pmatrix}e^\frac{t}{2}&0\\0&e^{-\frac{t}{2}}\end{pmatrix}$$ 
(for the appropriate normalization of the metric).

\begin{proof}[Proof of theorem \ref{sufficient conditions of density}] Let $g=(g_1,g_2)$ belong to $\rm G$ and denote by $z=(z_1,z_2)$ the basepoint in $\HH\times\HH$ of the Weyl chamber $g$.\\
\indent {\it 1)} We can assume (up to a permutation) that the fixed point is $g_1(\infty)$. This point is fixed by $h_1=p_1(\ga)$ where $\ga=(h_1,e_2)$ is a mixed element of $\GA$. The hyperbolic isometry $h_1$ of $\HH$ fixes the points $h_1^+$ and $h_1^-$ of $\dH$ and the elliptic isometry $e_2$ of $\HH$ fixes a point $x_2$ in $\HH$. Replacing, if necessary, $\ga$ by $\ga^{-1}$ we assume that $g_1(\infty)$ is the attractive fixed point $h_1^+$ of $h_1$. \\
We prove this first part in three steps, decreasing assumptions in each step.\\
\underline{Step 1:} We add here the following assumptions: $g_1(0)=h_1^-$ and $z_2=x_2$.\\
Since the element $\ga=(h_1,e_2)$ belongs to an irreducible lattice, the component $e_2$ is of infinite order. Therefore it generates a dense semi-subgroup of the stabilizer (conjugated to $\PSO(2,\RR)$) of the point $x_2$. If $k$ is any (elliptic) isometry of $\HH$ fixing the point $x_2$, there then exists a sequence $(m_n)_n$ of positive integers, going to $+\infty$, such that
$$\lim_{n\rightarrow+\infty}e_2^{-m_n}=k.$$
The point $g_1(\infty)$ (resp. $g_1(0)$) is the attractive (resp. repulsive) fixed point of the hyperbolic isometry $h_1$, therefore there exists an element $a_1$ in the diagonal semi-subgroup ${\rm A}_1^+$ of $\PSL(2,\RR)$ such that
$$h_1g_1=g_1a_1.$$
The element
$$(g_1,kg_2)=\lim_{n\rightarrow+\infty}\ga^{-m_n}g(a_1^{m_n},\Id)$$
belongs to the set $\ol{\GA g\rm A^+}$. Let $\tilde g_1$ be in $\PSL(2,\RR)$. By proposition \ref{refined density}, there exists a sequence $(\ga_n)_n$ in $\GA$ such that $(p_1(\ga_n))_n$ converges to $\tilde g_1g_1$ and $(p_2(\ga_n)z_2)_n$ converges to the point $g_2(0)$. For each $n$, there exists an elliptic oriented isometry $k_n$ of $\HH$ fixing $z_2$ such that $p_2(\ga_n)k_ng_2$ defines the unit tangent vector of $\HH$ based on the point $p_2(\ga_n)z_2$ and tangent to the geodesic line passing (in this order) through the points $p_2(\ga_n)z_2$ and $z_2$.
Therefore, for every $n$, the unit tangent vector
$$u_n=p_2(\ga_n)k_ng_2\phi^{d(z_2,p_2(\ga_n)z_2)}$$
is based on $z_2$ and defines an oriented geodesic line containing the point $p_2(\ga_n)z_2$. Since this last point converges to $g_2(0)$, the sequence $(u_n)_n$ converges to the unit tangent vector defined by $g_2$. We have
$$\lim_{n\rightarrow+\infty}\ga_n(g_1,k_ng_2)(\Id,\phi^{d(z_2,p_2(\ga_n)z_2)}) =(\tilde g_1,g_2)$$
which therefore belongs to $\ol{\GA g\rm A^+}$. We proved that $\PSL(2,\RR)\times\{g_2\}$ is contained in $\ol{\GA g\rm A^+}$. The fact that the second projection $p_2(\GA)$ is dense in ${\rm G}_2$ implies that $\ol{\GA g\rm A^+}$ equals $\rm G$, that is to say $x\rm A^+$ is dense in $\GA\backslash\rm G$.\\
\underline{Step 2:} Here we add only the assumption that $g_1(0)=h_1^-$.\\
By the same argument as in step 1, for every isometry $k$ fixing the point $x_2$, the element $(g_1,kg_2)$ belongs to the set $\ol{\GA g\rm A^+}$. By proposition \ref{refined density}, there exists a sequence $(\ga_n)_n$ in $\GA$ such that
$$\lim_{n\rightarrow+\infty}p_1(\ga_n)=\Id \quad\text{and}\quad \lim_{n\rightarrow+\infty}p_2(\ga_n)y_2=g_2(\infty)\;\text{for every point $y_2$ in $\HH$}.$$
The sequence $(p_2(\ga_n))_n$ is divergent in $\PSL(2,\RR)$, therefore the sequence $(p_2(\ga_n)^{-1}x_2)_n$ is divergent in $\HH$. For large enough $n$, the point $z_2$ is contained in the open hyperbolic disc of center $x_2$ and of radius $d(x_2,p_2(\ga_n)^{-1}x_2)$. Thus the geodesic ray (defined by $g_2$) based on $z_2$ and directed to $g_2(\infty)$ intersects the hyperbolic circle of center $x_2$ and passing through the point $p_2(\ga_n)^{-1}x_2$. This intersection is a point $k_np_2(\ga_n)^{-1}x_2$ where $k_n$ belongs to the stabilizer of $x_2$ in $\PSL(2,\RR)$ (Fig.1).
\vskip 5mm
\centerline{\hfill \includegraphics[height=50mm]{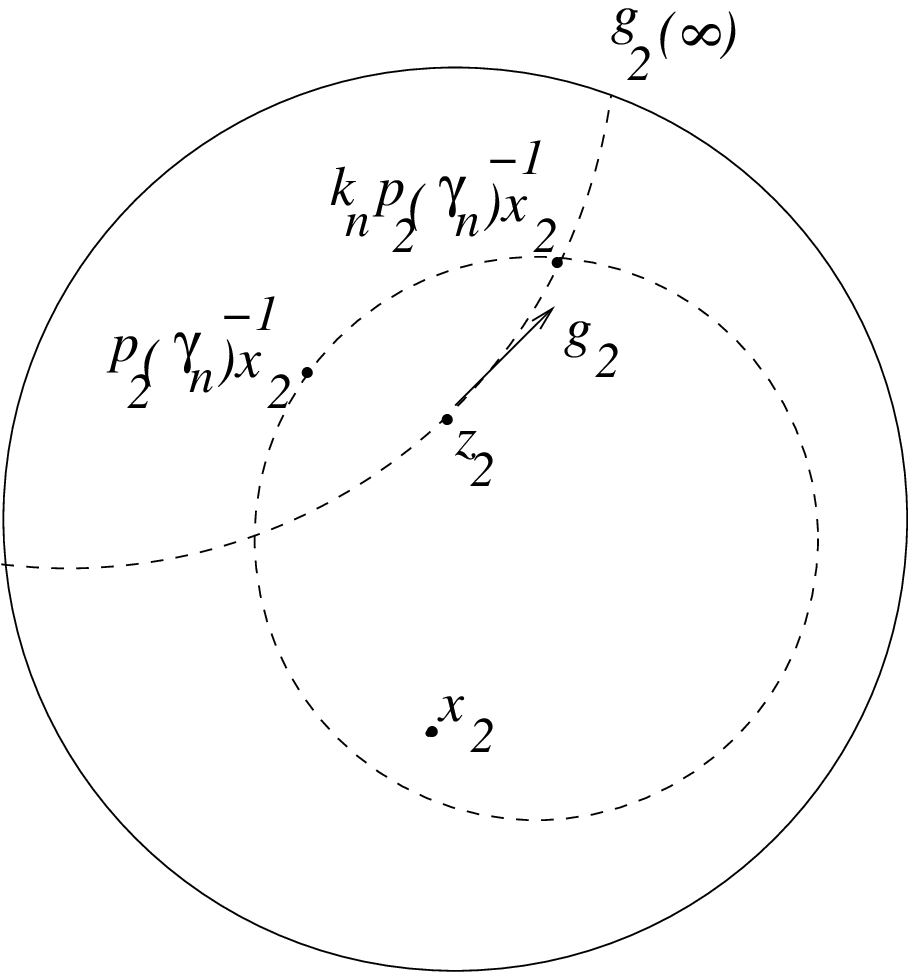}\hfill}
\leftline{\hfill Fig.~1\hfill}
\vskip 5mm
The geodesic ray defined by $p_2(\ga_n)k_n^{-1}g_2$ contains the point $x_2$. If $t_n$ equals the distance $d(x_2,p_2(\ga_n)k_n^{-1}z_2)$, then the element $p_2(\ga_n)k_n^{-1}g_2\phi^{t_n}$ defines for each $n$ a unit tangent vector based on the point $x_2$.
Passing to a subsequence, it converges to a unit tangent vector $\tilde g_2$ based on $x_2$ and satisfying
$$\lim_{n\rightarrow+\infty}\ga_n(g_1,k_n^{-1}g_2)(\Id,\phi^{t_n})=(g_1,\tilde g_2).$$
Therefore $\ol{\GA g\rm A^+}$ contains $\GA(g_1,\tilde g_2)\rm A^+$ which is dense in $\rm G$ by step 1.\\
\underline{Step 3:} No additional assumption, we only assume that $g_1(\infty)=h_1^+$.\\
Since the point $g_1(\infty)$ is distinct from the point $h_1^-$, $s=g_1^{-1}(h_1^-)$ is a real number and the element $u=\pm\begin{pmatrix}1&s\\0&1\end{pmatrix}$ of $\PSL(2,\RR)$ satisfies the following:
$$u^{-1}g_1^{-1}h_1g_1u=a_1=\pm\begin{pmatrix}\lambda&0\\0&\lambda^{-1}\end{pmatrix}$$
where $\lambda>1$ is the greatest eigenvalue of $h_1$. We have
$$u^{-1}g_1^{-1}h_1^{-n}g_1u=a_1^{-n},\quad h_1^{-n}g_1a_1^{n}=g_1ua_1^{-n}u^{-1}a_1^{n},$$
$$\text{and}\quad\lim_{n\rightarrow+\infty}h_1^{-n}g_1a_1^{n}=g_1u\quad\text{because}\quad \lim_{n\rightarrow+\infty}a_1^{-n}u^{-1}a_1^{n}=\Id.$$
Since the isometry $p_2(\ga)=e_2$ is elliptic, there exists a divergent sequence $(m_n)_n$ of positive integers such that
$$\lim_{n\rightarrow+\infty}e_2^{-m_n}=\Id.$$
We have
$$\lim_{n\rightarrow+\infty}\ga^{-m_n}g(a_1^{n},\Id)=(g_1u,g_2)$$
which is therefore an element of $\ol{\GA g\rm A^+}$. But this element $(g_1u,g_2)$ satisfies: $g_1u(\infty)$ and $g_1u(0)$ are fixed by $p_1(\ga)$. Therefore, by step 2, $\GA(g_1u,g_2)\rm A^+$ is dense in $\rm G$ and so is $\GA g\rm A^+$.\\

{\it 2)} We prove here the second part of the theorem. Now we assume that exactly one of the points $g_1(\infty)$ and $g_2(\infty)$ is fixed by a hyper-regular element $\ga=(h_1,h_2)$ of $\GA$. Let $h_i^\pm$ be the fixed points of $h_i$. We may assume that $g_1(\infty)$ is fixed by (the first component $h_1$ of) $\ga$ and $g_2(\infty)$ is not fixed by (the second component $h_2$ of) $\ga$.
The centralizer of $\ga$ in $\rm G$ is a conjugate of the subgroup ${\rm A}$; we denote it by ${\rm Z}={\rm Z}_1\times{\rm Z}_2$. Since the element $\ga$ of the irreducible lattice $\GA$ is hyper-regular, the discrete subgroup ${\rm Z}\cap\GA$ is a lattice in $\rm Z$ by theorem \ref{compact orbits}. The element $\ga$ is a positive power of a primitive element in the lattice ${\rm Z}\cap\GA$. Thus we can assume that $\ga$ is primitive: there exists another element $\ga'=(h'_1,h'_2)$ in the lattice ${\rm Z}\cap\GA$ such that ${\rm Z}\cap\GA$ is generated by $\ga$ and $\ga'$. For $i=1,2$, $h_i$ and $h'_i$ have the same fixed points $h_i^\pm$. Moreover we can assume they have the same attractive (resp. repulsive) point and their eigenvalues are such that the elements $\ga$ and $\ga'$ have the following position in (the Lie algebra of) $\rm Z$ (see Fig.2).
\vskip 5mm
\centerline{\includegraphics[height=50mm]{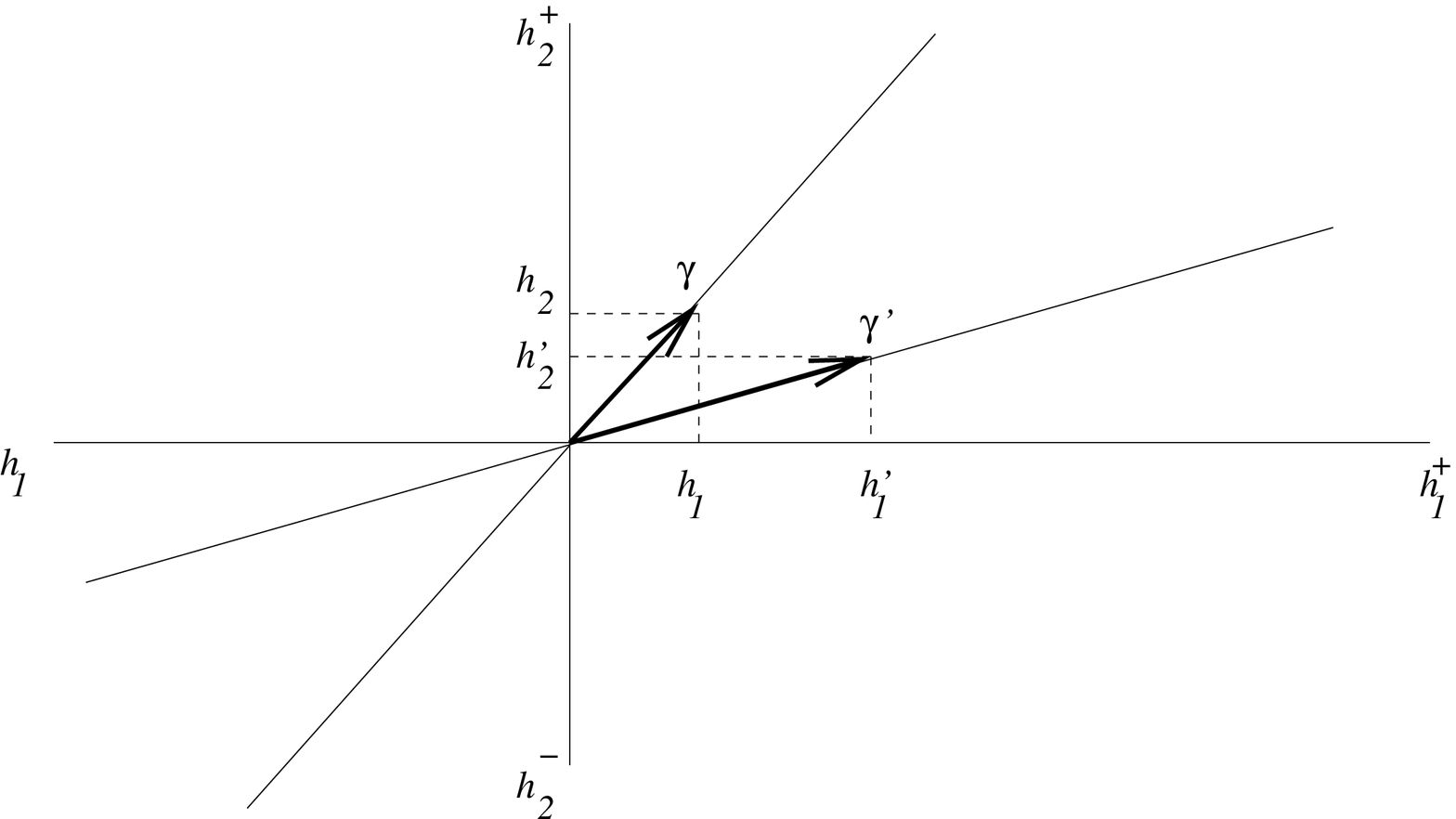}}
\leftline{\hfill Fig.~2\hfill}
\vskip 5mm

By irreducibility, the restriction to the lattice ${\rm Z}\cap\GA$ of each projection
$$p_i\;:\;{\rm Z}\longrightarrow{\rm Z}_i\quad\quad i=1,2,$$
is injective. Consequently the images $p_i({\rm Z}\cap\GA)$ are dense in ${\rm Z}_i$ and the closure of the sub-semigroup
$$\left\{h_2^m{h'_2}^{-m'}\;|\;m,m'\in\ZZ,m\geq0,m'\geq0\right\}$$
of the group ${\rm Z}_2$ has a non-empty interior. In the second factor $\dH$ of $\mathcal F$, the set of points fixed by a mixed element of the lattice $\GA$ is dense (it is non-empty and invariant under $p_2(\GA)$). Therefore there exists a mixed element $f=(f_1,f_2)$ in $\GA$, where $f_1$ is elliptic and $f_2$ is hyperbolic fixing a point $f_2^+$ which belongs to the set
$$\overline{\left\{h_2^m{h'_2}^{-m'}g_2(\infty)\;|\;m,m'\in\ZZ,m\geq0,m'\geq0\right\}},$$
(which has non-empty interior since $g_2(\infty)$ is fixed neither by $h_2$ nor by $h'_2$). There exist sequences $(m_n)_n$ and $(m'_n)_n$ of positive integers such that the sequence
$(h_2^{m_n}{h'_2}^{-m'_n})_n$ converges to an element $c$ of ${\rm Z}_2$ satisfying:
$$cg_2(\infty)=f_2^+.$$
By discreteness of $\GA$, the sequence of isometries $(h_1^{-m_n}{h'_1}^{m'_n})_n$ diverges (in ${\rm Z}_1$). The choice of $\ga$ and $\ga'$ implies moreover that, for $n$ large enough, the attractive fixed point of $(h_1^{-m_n}{h'_1}^{m'_n})_n$ is $h_1^+$. We now apply the same kind of argument than in step 3. The point $g_1(\infty)$ is distinct from the point $h_1^-$, thus $s=g_1^{-1}(h_1^-)$ is a real number and the element $u=\pm\begin{pmatrix}1&s\\0&1\end{pmatrix}$ of $\PSL(2,\RR)$ satisfies the following:
$$u^{-1}g_1^{-1}h_1^{-m_n}{h'_1}^{m'_n}g_1u= \pm\begin{pmatrix}\lambda_n&0\\0&\lambda_n^{-1}\end{pmatrix}=\phi^{2\ln\lambda_n}$$
where $(\lambda_n)_n$ is the sequence of the greatest eigenvalue of $h_1^{-m_n}{h'_1}^{m'_n}$. This sequence goes to $+\infty$ and we have
$$\lim_{n\rightarrow+\infty}(h_1^{-m_n}{h'_1}^{m'_n})^{-1}g_1\phi^{2\ln\lambda_n}=g_1u\quad \text{because}\quad\lim_{n\rightarrow+\infty}\phi^{-2\ln\lambda_n}u^{-1}\phi^{2\ln\lambda_n}=\Id.$$
Therefore we obtain
$$\lim_{n\rightarrow+\infty}\ga^{m_n}{\ga'}^{-m'_n}g(\phi^{2\ln\lambda_n},\Id)=(g_1u,cg_2)$$
and this limit belongs to the closure $\overline{\GA g\rm A^+}$. Since the point $cg_2(\infty)$ is fixed by the mixed element $f$ of $\GA$, then the point $(g_1u,cg_2)$ has a dense semi-orbit by $1)$. Therefore $\GA g\rm A^+$ is also dense in $\rm G$.
\end{proof}

\begin{remark} In step 1, we used the property that, if an elliptic isometry of $\HH$ fixing a point $z$ is of infinite order, it generates a dense subgroup in the stabilizer of $z$. This property is not true for elliptic isometries of higher-dimensional hyperbolic spaces. This is the obstruction to generalizing theorem \ref{sufficient conditions of density} to the Weyl chamber flow on irreducible quotients $\GA\backslash\HH^n\times\HH^n$ when $n\geq3$.
\end{remark}

\begin{proof}[Proof of the corollary \ref{corollary for uniform lattices}] Here we assume that $\GA$ is a uniform irreducible lattice and contains a non-trivial element $\ga$ which fixes one of the points $g_1(0)$,$g_1(\infty)$,$g_2(0)$ or $g_2(\infty)$. This element $\ga$ is semi-simple (because $\GA$ is uniform) but cannot be elliptic. Thus it is mixed or hyperbolic. If it is mixed, the orbit $x\rm A$ of the point $x=\GA g$ is dense by theorem \ref{sufficient conditions of density} and remark \ref{remarks in introduction} \textit{1)}. If it is hyperbolic, there exists an element $h=(h_1,h_2)$ in $\rm G$ such that $h^{-1}\ga h$ belongs to the subgroup $\rm A$. The orbit $y\rm A$ of the point $y=\GA h$ is closed by lemma \ref{properness}. Therefore it is compact and theorem \ref{compact orbits} implies the existence of a hyper-regular element $\ga'$ in $\GA$ which commutes with $\ga$. Thus $\ga'$ fixes the same points in $\cal F$. If all the points $g_1(0),g_1(\infty),g_2(0)$ and $g_2(\infty)$ are fixed by $\ga'$, the orbit is compact by theorem \ref{compact orbits}, otherwise the orbit is dense by theorem \ref{sufficient conditions of density} and remark \ref{remarks in introduction} \textit{1)}.
\end{proof}

\vskip 20mm
\noindent Damien \textsc{Ferte}\\
IRMAR -- Universit\'e Rennes I\\
35042 RENNES CEDEX\\
FRANCE\\
\textsl{damien.ferte@math.univ-rennes1.fr}

\end{document}